\documentclass[a4paper,10pt]{article}

\usepackage{amssymb,amsfonts,amsmath,geometry}

\newtheorem{proposition}{Proposition}

\newcommand{\cqd}{\nopagebreak\hfill\fbox{ }}

\usepackage{mathptmx}       
\usepackage{helvet}         
\usepackage{courier}        
\usepackage{type1cm}        
\usepackage{makeidx}         
\usepackage{graphicx}        
\usepackage{multicol}        
\usepackage[bottom]{footmisc}

\usepackage{amsfonts, amstext,amscd,amstext,latexsym}

\makeindex             

\begin{document}



\def\cqd {\,  \begin{footnotesize}$\square$\end{footnotesize}}
\def\cqdt {\hspace{5.6in} \begin{footnotesize}$\blacktriangleleft$\end{footnotesize}}
\def\refname{References}
\def\bibname{References}
\def\chaptername{\empty}
\def\figurename{Fig.}
\def\abstractname{Abstract}

\title{Transport and large deviations for Schrodinger operators and Mather measures}
\author{A. O. Lopes (*)  and Ph. Thieullen (**)}

%
\maketitle

\centerline{(*) \,Instituto de Matem\'atica-UFRGS, Avenida Bento Gon\c{c}alves 9500 Porto Alegre-RS Brazil}

\centerline{(**) Institut de Math\'ematiques, Universit\'e Bordeaux, F-33405 Talence, France.}

\abstract{In this mainly survey paper we consider the Lagrangian $ L(x,v) = \frac{1}{2} \, |v|^2 - V(x) $,
and a closed form $w$ on the torus $ \mathbb{T}^n $.
For the associated Hamiltonian we consider the
the Schrodinger operator ${\bf H}_\beta=\, -\,\frac{1}{2 \beta^2}   \, \Delta +V$ where $\beta$ is large real parameter. Moreover, for the given form $\beta\, w$ we consider the associated twist operator ${\bf H}_\beta^w$.
We denote by $({\bf H}_\beta^w)^*$ the corresponding backward operator.
We are interested in the positive eigenfunction  $ \psi_\beta$ associated to the
the eigenvalue
$ E_\beta$  for  the operator  ${\bf H}_\beta^{w} $. We denote $ \psi_\beta^*$ the positive eigenfunction associated to the
the eigenvalue
$ E_\beta$  for  the operator  $({\bf H}_\beta^{w})^* $.
Finally, we analyze the asymptotic limit of the probability $\nu_\beta=  \psi_\beta\, \psi_\beta^*$ on the torus when $\beta \to \infty$.
The limit probability is a Mather measure.
We consider Large deviations properties and we derive a result on Transport Theory.
We denote $L^{-}(x,v) = \frac{1}{2} \, |v|^2 - V(x) - w_x(v) $ and $L^{+}(x,v)
= \frac{1}{2} \, |v|^2 - V(x) + w_x(v) $. We are interest in the
transport problem from $\mu_{-}$ (the Mather measure for $L^{-}$) to
$\mu_{+}$ (the Mather measure for $L^{+}$) for some natural cost function. In the case the maximizing probability is unique we
use a Large Deviation Principle due to N.  Anantharaman in order to show that the conjugated sub-solutions $u$ and $u^*$ define an admissible  pair
which is optimal for the dual Kantorovich problem.}

\section{Introduction and basic definitions }\label{secaoinicial}

Given a closed form $w$ on the torus $ \mathbb{T}^n $ we consider the Lagrangian
$ L(x,v) = \frac{1}{2} \, |v|^2 - V(x) +w$, where $L:T \mathbb{T}^n \to \mathbb{R}$ and $T \mathbb{T}^n$ is the tangent bundle.

The infimum of $\int L(x,v) \, d\mu(x,v)$ among the invariant probabilities for the Euler Lagrange flow  on the
tangent bundle $T \mathbb{T}^n$  is called the critical value of $L$. A probability which attains such infimum is called a Mather measure (see \cite{CI} for references an general results).

We denote by $H(x,p)=\frac{1}{2} \, |p|^2 + V(x)  $ the associated Hamiltonian for  the Lagrangian $L(x,v)=\frac{1}{2} \, |v|^2 - V(x) $ and for each $\beta\in \mathbb{R}$  we consider the corresponding
Schrodinger operator ${\bf H}_\beta=\, -\,\frac{1}{2 \beta^2}   \, \Delta +
V$ for such Hamiltonian.

For each $\beta$ we consider a certain associated quantum state and quantum probability on $\mathcal{L}^2 (\mathbb{T}^n)$ (associated to an eigenvalue of ${\bf H}_\beta$) and we
are interested in the limit of such probability when $\beta \to \infty$.

We call $\beta$ the semiclassical parameter. In an alternative form we can take $\hbar=\frac{1}{\beta}$ and  consider the limit when $\hbar \to 0.$

An interesting relation of such limit probabilities with Mather measures was investigated by N. Anantharaman (see \cite{A3} \cite{A1} and \cite{A2})

We will present here  some of these results which are  related to transport and  large deviation properties.

\medskip

Consider $w(v) = <P,v>$ a closed form $w$ in the torus
$\mathbb{T}^n$, where $P$ is a vector in $\mathbb{R}^n$.

Suppose that $\mu_+$ and $\mu_{-}$ are respectively the Mather
measures for the Lagrangians
$$L^{+} (x,v) = \frac{1}{2}\, |v|^2 - V(x) + w(v)\,\, \mbox{and} \,\,L^{-} (x,v) = \frac{1}{2}\, |v|^2 - V(x) - w(v), $$

$ x\in \mathbb{T}^n$ and $V:\mathbb{T}^n \to \mathbb{R}$ smooth.

We assume the Mather measure is unique in each problem (see \cite{CI}, \cite{Fathi1}, \cite{Fathi}).

We will follow closely the notation of the nice exposition
\cite{A1} (see also \cite{A2} \cite{A3}). The results presented here in the future sections are inspired in \cite{LOP}. The main tool is the involution kernel introduced in \cite{BLT} (see  also
\cite{LM} \cite{LOS} \cite{LO} \cite{LMMS} \cite{LM2} \cite{CLO} \cite{LM} \cite{GLM} \cite{LR})

We will consider the Lax-Oleinik operator  $T^{-}_t$, $t\geq 0$,
given by
$$ T^{-}_t u_1(x) = \inf_{\gamma (t) =x, \gamma:[0,t] \to \mathbb{T}^n} \{ u_1(\gamma(0))  + \int_0^t \,L^{-}  (\gamma(s), \gamma '(s))ds  \}.$$

Denote by $u$ (Lipchitz), $ u : \mathbb{T}^n \to \mathbb{R}$, the
unique (up to additive constant because $\mu_{-}$ is unique)
solution of  $T^{-}_t u= u + t \,E$, for all $t\geq 0$, and where
$E$ is constant.

Consider the Lax-Oleinik operator  $T^{+}_t$, $t\geq 0$, given by
$$ T^{+}_t u_2(x) = \inf_{\gamma (0) =x, \gamma:[0,t] \to \mathbb{T}^n} \{ u_2(\gamma(t))  - \int_0^t \,L^{-}  (\gamma(s), \gamma '(s))ds  \}.$$

Denote by $u^*$ the Lipchitz  function $ u^* : \mathbb{T}^n \to
\mathbb{R}$,  such that,  $T^{+}_t (-u^*)= -u^* + E \, t .$

We assume that $u$ and $u^*$  are such $u+ u^*$ is zero in the
support of $\mu_{-}$.

The function $W(x,y)$ (which could be called the convolution kernel) is
given by the bellow expression

$$-\,\inf_{\, \alpha\,\, \in \,C^1([0,1], \mathbb{T}^n), \, \alpha
(0) = x ,\, \alpha(1) =y\,}\{\int_0^1 \,\, [\,-\,V(\alpha(s) +
w(\alpha '(s))\, ]\, ds \,+\, \int_0^1 \frac{1}{2}\,
||\alpha'(t)||^2\}.$$

We denote by $h(y,y)$ the Peierls barrier for the Lagrangian
$L^{-}$. In the present case $ h(y,y) = u(y) + u^* (y).$

Main references on Transport Theory are \cite{Vi1} \cite{Vi2} and \cite{Ra}.

We denote by ${\cal K}( \mu_+, \mu_{-})$ the set of  probabilities
$\hat{\mu}$ on $ \mathbb{T}^n \times  \mathbb{T}^n$, such that
respectively $ \mu_+ =\pi_1^\# ( \hat{\mu}    ) $ and $ \mu_{-}=
\pi_2^\# ( \hat{\mu}    ).$

Give $c(x,y)$ we say that $f$ and $g$ are $c$-admissible if, for any
$x,y \in \mathbb{R}^n$, we have $ f(x) - g(y) \leq c(x,y).$ We
denote by ${\cal F}$ the set of such pairs $(f,g)$.

We will consider, for  the cost function $c(x,y) = -\, W(y,x) $, a
$c$-Kantorovich problem

$$ \inf_{  \hat{\mu}  \in  {\cal K}( \mu_+, \mu_{-}) }\, \int \int
c(x,y)\,d\hat{\mu} (x,y).$$

We denote the minimizing probability by $\hat{\mu}_{\min}$. Note
that this probability  projects on the second variable on $
\mu_{-}$.

Note that the transport optimal probability for $-W$ and for $-W+I$ (where $I$ is the Peierl's barrier) are the same.

We point out that the projected projected Mather measures $\mu^+$ and $\mu^{-}$  are the same in the present case.

We will show here that the dual problem for $-W$
$$ \max \{\, \,\int  \,f(x) \, d\mu_+ (x)\,-\, \int\,g(y) \, d\mu_{-}(y)\,\,|\,\, f(x) - g(y) \leq c(x,y)\,
\}=$$
$$ \max \{\, \,\int  \,f(x) \, d\mu_+ (x)\,-\, \int\,g(y) \, d\mu_{-}(y)\,\,|\,\, (f,g) \in  {\cal F}\,
\},$$

has a pair of optimal solutions $(u, u^*)$ which are the viscosity
solutions of the Hamilton-Jacobi equations (fixed points of the
corresponding Lax-Oleinik operators as defined above)

We can consider alternatively (the same problem)

$$ \inf_{  \hat{\mu}  \in  {\cal K}( \mu_+, \mu_{-}) }\, \int \int
\tilde{c}(x,y)\,d\hat{\mu} (x,y),$$

where $\tilde{c} (x,y) = -W(y,x) + h(y,y)$. The introduction of a
function on the variable $y$ which vanishes in the support of
$\mu_{-}$ does not change the minimizing measure. However, this new problem have a different optimal pair.

We denote by $ {\cal W}_ {x}^h  $ the Brownian motion in
$\mathbb{R}^n$ (with coefficient $h$, that is, at time $t=1$ the
variance is $\sqrt{h}$) beginning at $x$, and $ {\cal W}_ {x,y,t}^h$
its disintegration at the point $y$ and at the time $t$.

Consider the Schrodinger ${\bf H}^h=\, -\, \frac{h^2}{2} \, \Delta +
V$ (where $V$ is the periodic extension to $\mathbb{R}^n$) which
acts on real (periodic) functions defined in $\mathbb{R}^n$. It is
known that ${\bf H}$ has pure point spectrum (see \cite{Dav} and \cite{LS}).

Note that
$$-\,\frac{1}{h} {\bf H}^h\,=\, \frac{h}{2} \, \Delta - \frac{1}{h}
V.$$

The Kernel $K(x,y,t)$ of the extension of $e^{-\,\frac{t}{h} \,H}$
to an integral operator is (see \cite{A1})
$$ K(x,y,t) \,=\, \int e^{  - \,\frac{1}{h} \, \int_0^t \, V( \alpha (s))
\,ds}\,  {\cal W}_ {x,y,t}^h (d\, \alpha).$$

Given
$$L^{w} (x,v) = \frac{1}{2}\, v^2 - V(x) \,- \,w(v)  =   \frac{1}{2}\, v^2 - V(x) \, -\, <P,v> ,$$ the
corresponding Hamiltonian $H^w(x,p)$ via Legendre transform is
$$ H^w (x,p) = \frac{|| p\, +\, P\,||^2}{2}\, + \, V(x).$$

In the same way, for
$$L^{+} (x,v) = \frac{1}{2}\, v^2 - V(x) \,+ \,w(v)  =   \frac{1}{2}\, v^2 - V(x)\, -\, <P,v> ,$$ the
corresponding Hamiltonian $H^{w^*}(x,p)$ is
$$ H^{w^*} (x,p) = \frac{|| p\, -\, P\,||^2}{2}\, + \, V(x).$$

Consider, a certain point $x_0={\cal O}\in \mathbb{R}^n$ fixed (on
the universal cover of the torus). As the form $w$ on the torus is
closed, it is exact on the lifting to the universal cover, then, the
value $ \int^{x}_{x_0} \,w $ does not depend on the path we choose
to connect $x_0$ to $x$.
\bigskip

\section{Transport in the configuration space for the
Aubry-Mather problem}

For each real value $\beta$ we consider the operator $$  {\bf H}_
\beta^w \,=\, e^{- \beta \, \int^{x}_{x_0}\, w   }\, \circ \, {\bf
H}_\beta \, \circ e^{\beta\, \int^{x}_{x_0} \, w }  \,=\, e^{- \beta
\, \int^{x}_{x_0}\, w }\, \circ (\, -\,\frac{1}{2 \beta^2} \, \Delta
\,+\, V\,)\, \circ e^{\beta\, \int^{x}_{x_0} \, w } .$$

We can consider such operator acting on the torus or on
$\mathbb{R}^n$. When we consider the Brownian motion we
should consider, off course, its action on $\mathbb{R}^n$.

The Kernel $K(x,y,t)$ of the extension of $e^{\,t \beta \,{\bf
H}_\beta^w}$ to an integral operator is given by
$$ K_\beta (x,y,t) \,=\, \int e^{\, -\, \beta \, \int_0^t \, \,V( \alpha
(s))\,ds\, + \beta\, <\,P \,,\, (y-x)\,>\,}\,  {\cal W}_
{\,x,y,t}^{\,\beta^{-1}} (d\, \alpha).$$

Note that above we consider the integral
$$ \beta \,\int_0^t \,
[\,-\,V( \alpha (s))\, + \,w_{  \alpha (s) }\,(\alpha' (s)) \,]\,
ds.$$

${\bf H}_\beta^w $ is not self adjoint but has a real pure point
spectrum.

We denote by $ E_\beta$ the maximum eigenvalue of    ${\bf
H}_\beta^w $ (acting on real functions) and $ \psi_\beta$ is the corresponding normalized real
eigenfunction in ${\cal L}^2 ( \mathbb{T}^n, dx)$. The positive eigenfunction
$ \psi_\beta$ is unique if we assume its norm is $1$. It's the only
totally positive eigenfunction of ${\bf H}_\beta^w $ (see \cite{A1}
expression (3.15)). The eigenvalue is simple and isolated (see appendix on \cite{A2}).

For each real value $\beta$ we consider the $w$-backward operator $$
{\bf H}_\beta^{w^*} \,=\, e^{\beta \, \int^{x}_{x_0}\, w   }\, \circ
\, {\bf H}_\beta \, \circ e^{-\beta\, \int^{x}_{x_0} \, w }\,=\, e^{
\beta \, \int^{x}_{x_0}\, w }\, \circ (\, -\,\frac{1}{2 \beta^2} \,
\Delta \,+\, V\,)\, \circ e^{-\beta\, \int^{x}_{x_0} \, w } .$$

We will be interested in high values of $\beta$.

$ E_\beta$ is the maximum eigenvalue of    ${\bf H}_\beta^{w^*} $
and we denote by $ \psi_\beta^* $ the corresponding real
eigenfunction in ${\cal L}^2 ( \mathbb{T}^n, dx)$. Similar
properties to the case of  $\psi_\beta$ are true for such
eigenfunction. We assume   $\int\, \psi_\beta^* (x) \, dx=1$ and
also $ \int \psi_\beta(x) \, \psi_\beta^* (x) \, dx=1$.

${\bf H}_\beta^{w^*} \, \circ \, {\bf H}_\beta^w $ is self adjoint.

We will be interested here in the probabilities
$$\nu_\beta (dx) \,=\,\psi_\beta (x)\, \psi_\beta^*(x) \, dx.$$

The probability $\nu_\beta (dx) \,=\,\psi_\beta (x)\,
\psi_\beta^*(x) \, dx$ is  stationary for the Markov operator
$$Q^t (f) (x) = e^{-\, t\, E_w} \,\psi_\beta (x)^{-1} e^{ t \,  {\bf H}_\beta^w}  (\psi_\beta \, f) (x)  $$

on the torus $\mathbb{T}^n$ (see \cite{A2}).

The correct point of view is to consider
$\psi_\beta$ as an eigenfunction and  $ \rho_\beta \,=\,
\psi_\beta^*(x) \, dx $ as an eigen-probability  for the semi-group
$t \, \to \,e^{ t \, {\bf H}_\beta^w}$.

Consider
$$ u_\beta = - \frac{ \log  \psi_\beta  }{\beta   }\,\,\,\,
\text{and}\,\,\,\, u_\beta^*  = - \frac{ \log  \psi_\beta^*  }{\beta
}.$$

It is known that the following equalities are true:

$$ - \, \frac{1}{2\, \beta} \, \Delta \, u_\beta \, +\, H^w (x,
d_x u_\beta)= E_\beta,$$ and
$$ - \, \frac{1}{2\, \beta} \, \Delta \, u_\beta^* \, +\, H^w (x,-
d_x u_\beta^*)= E_\beta,$$

The $\beta$-families of functions  $u_\beta$ and $u_\beta^*$ are
equi-Lipschizians and we can obtain from this fact convergent
subsequences. We assume here  the Mather measure is unique, and
therefore the limits exist in the uniform convergence topology, that
is
$$ \lim_{\beta \to \infty} \, u_\beta \,=u \, \, \,\,\text{and}\,\,\,\, \lim_{\beta \to \infty} \,
u_\beta^* \,=u^*.$$

It is known that $\lim_{\beta \to \infty}\, E_\beta$ exist and we
denote this value by $E$.

By stability of the viscosity solutions, the limits $u$ and $u^*$ are,
respectively,
viscosity solutions of the equations
$$H^w ( x, d_x \,u) = E\,\,\,\,\text{and}\,\,\,\, H^w ( x, -d_x \,u^*) = E $$

We assume also that the Mather measure $\mu$ for the lagrangian
$L^w$ is unique. In this case it is known (see for instance
\cite{A1} \cite{A2}) that in the weak topology
$$\lim_{\beta \to \infty}\, \nu_\beta\,=\, \mu.$$

In proposition 3.11 in \cite{A1} the following Large Deviation
Principle is obtained (see also \cite{A2} \cite{A3}):

\begin{proposition}
 Suppose the Mather measure is unique. Suppose also that in the uniform convergence topology
$$ \lim_{\beta \to \infty} \, u_\beta \,=u \, \, \,\,\text{and}\,\,\,\, \lim_{\beta \to \infty} \,
u_\beta^* \,=u^*.$$

Then, for $I(x) = u(x) + u^* (x)$ (from the normalization we choose
before $I(x)\geq 0$), we have

1) for any open set $O\subset \mathbb{T}^n$,

$$ \liminf_{\beta \to \infty} \, \frac{1}{\beta}\, \log \, \nu_\beta (O)\,=-\, \inf_{x\in
O}\,\{ I(x)\},$$

and,

2) for any closed set $F\subset \mathbb{T}^n$,

$$ \limsup_{\beta \to \infty} \,  \frac{1}{\beta}\, \log \, \nu_\beta (F)\,=\, -\inf_{x\in
F}\,\{ I(x)\}.$$
\end{proposition}

It follows from Varadhan's  Integral Lemma (section 4.3 in
\cite{DZ}) that, for any $C^\infty$ function  $F(x)$,
$$ \lim_{\beta \to \infty}\,  \frac{1}{\beta}\, \log \,\int \, e^{\beta\, F(x)}\, d \, \nu_\beta (x)\, =\, \sup_{x\in \mathbb{T}^n}
\,\{F(x) - I(x)\} .$$

The $W_\beta^t$-Kernel is defined by

$$ e^{ W_\beta^t (y,x)} =  \int e^{\,-\,   \beta \, \int_0^t \, \,V( \alpha
(s))\,ds\, - \beta \,   \int_y^x w  \,}\,  {\cal W}_
{\,y,x,t}^{\beta^{-1}} (d\, \alpha)= $$ $$ \int e^{\,  -\, \beta \,
\int_0^t \, \,V( \alpha (s))\, d s \, - \beta \, <\,P \,,\,
(x-y)\,>\,}\,  {\cal W}_ {\,y,x,t}^{\beta^{-1}} (d\, \alpha).$$

Note the plus sign on $V$.

Note that we exchange $x$ and $y$ above (with respect to the
previous considerations).

It is known (see \cite{A1}) that for any $\beta$ and any $t$
$$  \psi_\beta  (x) = \int\,\,e^{ W_\beta^t (y,x)}\, \psi_\beta^* (y) \, dy=
\int\,\,e^{ W_\beta^t (y,x)}\, \frac{1}{\psi_\beta (y)} d \nu_\beta
(y)
$$

Now from Schilder's Theorem and Varadhan's Integral Lema (see
\cite{DZ} also Theorem 4.3.9 in \cite{A2})
$$ - \,W(y,x)\,:=- \, \lim_{\beta\to \infty}  \, \frac{1}{\beta} \, \log e^{ W_\beta^{\frac{1}{\beta}} \,(y,x)}
=$$
$$  \inf_{\, \alpha\,\, \in \,C^1([0,1], \mathbb{T}^n),
\, \alpha (0) = y ,\, \alpha(1) =x\,}\{\int_0^1 \,\,
[\,-\,V(\alpha(s) + w_{\alpha (s)} (\alpha '(s))\, ]\, ds \,+\,
\int_0^1 \frac{1}{2}\, ||\alpha'(t)||^2\}
$$

Note above the plus sign on $w$.

The function $W(y,x)$ is the function $-I(y,x)$ in the notation of
\cite{A2}.

For any $\beta$
$$\frac{1}{\beta} \log \,(  \psi_\beta (x))\, =\,\frac{1}{\beta} \log \,( \int\,
e^{ W_\beta^{\frac{1}{\beta}}\, (y,x)}\,  \psi_\beta (y)^{-1}   \,d
\nu_\beta (y)).$$

Taking limits as $\beta\to \infty$ and using the Varadhan's integral
Lemma once more we get
$$ -u (x) = \sup_{z \in \mathbb{T}^n}\{    W(z,x) + u (z)            - I(z)\}.$$

Therefore, for any $x,y$ we get that
$$ -u(y) - u(x) \,\geq \, W(y,x)             - I(y).$$

From this we get
\begin{proposition}
$$ u(y) + u(x) \,\leq\,  -W(y,x)             + I(y) \,=\,c(y,x).$$

and the pair $ (u, u)$ is (-$W$+I)-admissible.
\end{proposition}

In the same way
the pair $ (u, u^*)$ is -$W$-admissible.

\begin{proposition} If $\hat{\eta}$ is an optimal minimizing
transport probability for $c$ and if $(f,f^\#)$ is an optimal pair
in ${\cal F}$, then the support of $\hat{\eta}$ is contained in the
set
$$\{\, (x,y) \in M \times M \, \,\mbox{such that}\,\, (f(x) - f^\#
(y))\, = \, c(x,y)\,\}.$$
\end{proposition}

{\bf Proof:} It follows from the primal and dual linear programming
problem formulation. The condition above is called the complementary
slackness condition (see \cite{EG}). \cqd

If one finds $\hat{\eta}$ an an admissible pair  $(f,f^\#)$
satisfying the above claim (for the support) one solves the
Kantorovich problem, that is, one finds the optimal transport
probability $\hat{\eta}$ .

From the above it follows.
\begin{proposition}
{\bf Proposition:} For $(x,y)$ in the support of $\hat{\mu}$ we have

$$u(x) + u(y) \,= \,  -W(y,x)             + h(y,y)= c(x,y),$$

or

$$u(x) + u (y) \,= \,  -W(y,x)  \, + \, ( u^* (y) +      u(y)) .$$

This means, for $(x,y)$ in the support of $\hat{\mu}$

$$u(x) - u^* (y) \,= \,  -W(y,x)  \,  .$$
\end{proposition}

 In other words, for any $x,y$   in the support of $\hat{\mu}_{min}$ we have that
 $ u(x)$ is given by
 $$  \inf_{\, \alpha\,\, \in
\,C^1([0,1], \mathbb{T}^n), \, \alpha (0) = y ,\, \alpha(1)
=x\,}\{\int_0^1 \,\, [\,-\,V(\alpha)+ w(\alpha ')\, ]\, ds
\,+\, \int_0^1 \frac{1}{2}\, ||\alpha'||^2ds\}\,+\,   u^* (y).$$

\end{document}